\documentclass{elsarticle}
\usepackage{amsfonts}
\usepackage{amsmath}

\newtheorem{theorem}{Theorem}

\begin{document}

\title{Coalescent approximation for structured populations in a stationary random environment}
\author[1]{S.Sagitov}
\author[1]{P.Jagers}
\author[2]{V.Vatutin}
\address[1]  {Mathematical Statistics, Chalmers University of Technology and University of Gothenburg, SE-412 96 Gothenburg, Sweden.}
\address[2]{Steklov Institute of Mathematics, Moscow, Russia}
\begin{abstract}
We establish convergence to the Kingman coalescent for the genealogy of a geographically - or otherwise - structured version of the Wright-Fisher population
model with fast migration. The new feature  is that migration
probabilities may change in a random fashion. This brings a novel formula for the coalescent effective population size (EPS). We call it a quenched EPS to emphasize the key feature of our model -  random environment. The quenched EPS is compared with an annealed (mean-field) EPS which describes the case of constant migration probabilities obtained by averaging the random migration probabilities over possible environments.
\end{abstract}
\maketitle

\section{Introduction}
The Wright-Fisher population model is used as a benchmark to measure the speed of the random genetic drift in actual biological populations as well as in population models with more structure than the classical setup allows  \cite{E}. Viewed backward in time, it is approximated by the Kingman coalescent, a simple algorithm of consecutively joining together pairs of sampled ancestral lines until a random ancestral tree is formed. The resulting process \cite{K1} has no parameters and the Wright-Fisher population size $N$ is mirrored in the time scale ensuring the coalescent approximation. The larger is $N$, the slower the rate of genetic drift, since it takes longer for an allele to get fixed in the population - in the coalescent tree this is reflected in longer branch lengths (as counted in generations).

If the genealogy of another, usually more structured, population model is approximated by the standard Kingman coalescent, then the time scale of the latter  $N_e$ takes the role of the Wright-Fisher population size. This is why it is called the {\it coalescent effective population size} (see \cite{SKK} as well as \cite{M2} and \cite{NK}). The effective size $N_e$ is usually smaller than the actual population size $N$ as $N_e$ incorporates a number of factors not present in the Wright-Fisher model that increase variability in the underlying genetic sampling process and thereby speed up genetic drift. Such factors might be demographic fluctuations \cite{JS} or age-structure \cite{SJ}. The recent note \cite{WS} discusses extensions of the coalescent effective population size concept. 

In settings where no coalescent approximation avails itself, ideas become more complicated, and several definitions circulate in literature (see \cite{E82} and \cite{E89}). Among these, the so called inbreeding effective population size (Crow and Kimura \cite[p.~347]{crowkimura} and Ewens \cite{E79}) is the one  that is closest in spirit to the coalescent effective population size.

A case studied by several authors (see \cite{Na}, \cite{No}) and nicely summarized in \cite{NK} is that of a geographically structured Wright-Fisher model with fast migration. It deals with a population living on $L\ge2$ islands with a constant total population size $N$ and where also population sizes on the islands  $Na_1,\ldots,Na_L$ are constant over time. The fixed population structure is then described by the positive vector \begin{equation}\label{a}
(a_1,\ldots,a_L),\ a_1>0,\ldots,a_L>0,\ a_1+\cdots+a_L=1.
\end{equation}

Let $b_{ij}$ denote the probability that a lineage located on island $i$ comes from island $j$ if traced one generation back in time. Clearly $\sum_{j=1}^Lb_{ij}=1$. If the backward migration matrix  $\mathbf{B}_1=(b_{ij})$ has a stationary distribution $(\gamma_1,\ldots,\gamma_L)$, 
%with \[\gamma_1\ge0,\ldots, \gamma_L\ge0,\ \gamma_1+\cdots+\gamma_L=1,\]
the ancestral process converges (see Section 2.2 in \cite{NK}) to the Kingman coalescent, provided time is scaled by the factor $N_e=N/c_f$, where
\begin{equation}\label{cf}
  c_f=\sum_{k=1}^L{1\over a_k}\gamma_k^2.
\end{equation}
It is easy to interpret the factor $c_f$ in $N_e =N/c_f$: two lineages coalesce, if while visiting the same island $k$ they both chose the same parent among $Na_k$ available.

In cases of slow migration (when the ancestral process is approximated by the {\it structured coalescent}) the effective population size formulae may give the impression that the effective population size significantly exceeds actual size (\cite{NT} and \cite{W98}). This phenomenon can be viewed as an artifact of the random sampling design: if two lineages are sampled from different sub-populations, it takes some time before they enter the same sub-population and get a chance to merge.

\begin{figure}\label{motex}
\centering
\includegraphics[width=10cm]{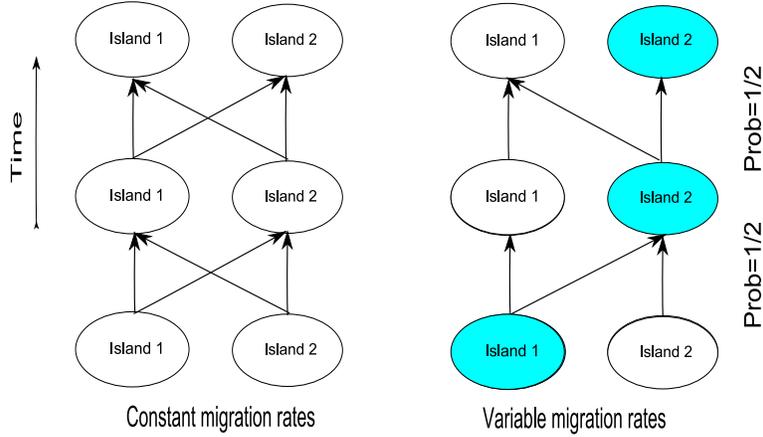}
\caption{Two-island modifications of the Wright-Fisher model.}
\end{figure}

We take a further step towards  more realistic models by allowing {\it variable migration} probabilities. The idea is illustrated in Figure \ref{motex}, presenting two versions of two-island populations (i. e. $L=2$). The right panel depicts a situation where for a given year each of the two islands can have an environmental advantage with equal probabilities, the advantage being that the offspring from the favored part can migrate to the other island but not vice verse. The left panel represents the corresponding constant environment case obtained by averaging over environmental fluctuations.

Our main result, Theorem \ref{main} in Section \ref{Sv}, on geographically structured populations with variable migration can be summarized as follows. If the backward migration matrix $\mathbf{B}_1$ is random, then the stationary distribution $(\gamma_1,\ldots,\gamma_L)$  also becomes a {\it random vector} and the coalescent effective population size formula takes the form
\begin{equation}\label{cq}
N_e=N/c_q,\  c_q=\sum_{k=1}^L{1\over a_k}\mathbb{E}(\gamma_k^2).
\end{equation}
Here the expectation operator is taken with respect to the randomly varying environment. If the random stationary probabilities are directly averaged into $\bar{\gamma}_k = \mathbb{E}(\gamma_k)$, and then inserted in (\ref{cf}), the result is an {\it annealed} (or in physics language {\it mean-field}) expression,
\begin{equation}\label{ca}
  c_a=\sum_{k=1}^L{1\over a_k}\bar{\gamma}_k^2.
\end{equation}

%(Warning: we use the same notation, most notably $(\gamma_1,\ldots,\gamma_L)$, for the random and constant environment cases. This makes clearer the parallel between the two cases but also requires some extra attention while reading the results.)

The expressions $c_a$ and $c_q$ thus pertain to the annealed and {\it quenched} approaches, respectively. Formula \eqref{ca} is interpreted as applied to the population with a constant environment obtained by averaging over all possible environmental scenarios (left panel in Figure \ref{motex}). The difference between  $c_a$ and $c_q$ is given by a weighted sum of variances
\begin{equation}\label{cqa}
c_q-c_a=\sum_{k=1}^L{1\over a_k}Var(\gamma_k).
\end{equation}
Formula (\ref{ca}) and Jensen's inequality imply $c_a\ge1$:
$$\sum_{k=1}^L{1\over a_k}\bar\gamma_k^2=\sum_{k=1}^La_k\left({\bar\gamma_k\over a_k}\right)^2\ge\left(\sum_{k=1}^La_k{\bar\gamma_k\over a_k}\right)^2=1.$$
This observation together with (\ref{cqa}) yields the important inequalities
\[1\le c_a\le c_q,\]
saying that
\[N^{\rm [quenched]}_e\le N^{\rm [annealed]}_e\le N.\]
According to (\ref{cqa}) the quenched and annealed effective population sizes coincide, $N^{\rm [quenched]}_e=N^{\rm [annealed]}_e$, if and only if the environment is constant, so that all $Var(\gamma_k)=0$. The quenched $N_e$ becomes strictly smaller than the annealed $N_e$, if there is an extra source of variability in genetic sampling due to random environment. Observe also that the effective population size is equal to the actual size $N$ only if migration probabilities faithfully follow the given population structure in that $\gamma_k=a_k$ for all $k=1,\ldots,L$. This holds, for example, in the ``dummy island'' case corresponding to the standard Wright-Fisher model (discussed as a test example in Section \ref{Sv}).

After this overview, the paper is organized as follows. Section \ref{Sv} contains a full description of the population model in a stationary random environment and the main result of the paper, Theorem \ref{main}, on convergence to the Kingman coalescent. Section \ref{ex} presents two detailed examples illustrating Theorem \ref{main} in the case of iid random environment. In Section \ref{Sc} we outline the main idea of the proof of the annealed $N_e$-factor formula \eqref{ca} given in \cite{NK}, using terms to which we shall refer in our analysis of  variable migration in Section \ref{pro}.\\

\textbf{Acknowledgements.} {\it This work has been supported by The Bank of Sweden Tercentenary Foundation and the Centre for Theoretical Biology at the University of Gothenburg. The third author was supported in part by the Russian Foundation for Basic Research, grant 08-01-00078. We thank three anonymous referees whose suggestions helped to improve the presentation significantly.}

\section{Convergence to the Kingman coalescent}\label{Sv}

The standard Wright-Fisher model with a constant population size $N$ represents an idealized population, lacking any kind of structure. The Wright-Fisher reproduction rule says that $N$ children are allocated to $N$ available parents uniformly at random.  Let $X(u)$ be the number of ancestral lineages $u$ generations backwards in time when $X(0)=n$ individuals were randomly sampled from the Wright-Fisher population. The time homogeneous Markov chain $\{X(u)\}$ with the finite state space $\{1,\ldots,n\}$ has a transition matrix $\mathbf{\Pi}=\mathbf{\Pi}_N$ such that
\begin{equation}
  \label{k}
\mathbf{\Pi}=\mathbf{I}+N^{-1}\mathbf{Q}+\mathbf{o}(N^{-1}),\ N\to\infty.
\end{equation}
Here $\mathbf{I}$ is the unit matrix of appropriate size, $\mathbf{o}(N^{-1})$ stands for a matrix whose elements are all of size $o(N^{-1})$, and $\mathbf{Q}=(q_{ij})_{i,j=1}^n$ with
\begin{equation}
  \label{Q}
 q_{ii}=-{i\choose2},\ q_{i,i-1}={i\choose2},
\end{equation}
and $q_{ij}=0$ whenever $i\ge j+2$ or  $j\ge i+1$. Thus  ($[x]$ standing for the integer part of $x$),
\begin{equation}
  \label{11}
\mathbf{\Pi}^{[Nt]}\to e^{t\mathbf{Q}},\ N\to \infty,
\end{equation}
implying the weak convergence (see the remark in the end of this section)
\begin{equation}
  \label{c1}
\{X([Nt]), t\ge0\}\to\{K(t), t\ge0\},\ N\to \infty,
\end{equation}
to a pure death process $\{K(t)\}$ with the infinitesimal transition matrix $\mathbf{Q}$. In view of \eqref{Q}, the latter means that $K(t)$ stays at the current state $i$ for an exponential time with mean $1/{i\choose2}$ and then jumps to $i-1$, until it is absorbed at $i=1$. This is the essence of the Kingman coalescent approximation for the standard Wright-Fisher model \cite{K1}.

As mentioned in the introduction, an important modification of the Wright-Fisher model adds a geographical structure, dividing the population of size $N$ into $L\ge2$ sub-populations of constant sizes $Na_1,\ldots,Na_L$, $a_1+\cdots+a_L=1$.
%(living on different islands). Assume (cf \cite{NK}) that the sub-population 
%are also constant over time so that the fixed population structure is described by a positive vector $(a_1,\ldots,a_L)$ satisfying  
Suppose, a lineage located on island $i$ may lead to island $j$, if followed one
generation back in time, with probability, say $b_{ij}$.
If the backward migration matrix $\mathbf{B}_1=(b_{ij})$ 
%with $$\sum_{j=1}^Lb_{ij}=1,\ i=1,\ldots,L$$
has a stationary distribution $(\gamma_1,\ldots,\gamma_L)$, then it is known (see Section \ref{Sc}) that
\begin{equation}
  \label{ce}
\{X([Nt/c_f]), t\ge0\}\to\{K(t), t\ge0\},\ N\to \infty,
\end{equation}
where  $c_f$ is defined by \eqref{cf}.

As a test case, consider again the standard Wright-Fisher model with $N$ individuals labeled by $1,\dots,N$ in any given generation. For a given vector \eqref{a} introduce a dummy island structure by assigning individuals
$$[N(a_1+\cdots+a_{i-1})]+1,\ldots,[N(a_1+\cdots+a_i)]$$
to the $i-$th island, $i=1,\ldots,L$, where $a_1+a_0=0$. Notice that in this case the backward migration probabilities depend on $N$ in the following weak way
\begin{equation}
  \label{N}
\mathbf{B}_1(N)=\mathbf{B}_1+N^{-1}\mathbf{D}_1(N).
\end{equation}
Here the main term matrix
\begin{equation*}
 \mathbf{B}_1=
\left(
\begin{array}{ccc}
a_1&\ldots&a_L\\
\vdots&\ldots&\vdots\\
a_1&\ldots&a_L
\end{array}\right)
\end{equation*}
readily gives the stationary distribution. The discrepancy matrix $\mathbf{D}_1(N)$ has negligible effect (see Appendix B), since the absolute values of its elements
\[d_{ij}=[N(a_1+\cdots+a_i)]-[N(a_1+\cdots+a_{i-1})]-a_iN\]
are all bounded by a constant independent of $N$. The insertion $\gamma_i=a_i$ into \eqref{cf} gives $c_f=1$, as it should.

We render the previous model more flexible by allowing the migration probabilities to change randomly from generation to generation. Let $b_{ij}^{(u)}$ denote the probability that a
lineage located on island $i$ at the backward time $u-1$ comes from island $j$, if followed one further generation
back in time, so that $\sum_{j=1}^Lb_{ij}^{(u)}=1$. We will treat the backward migration matrix
$\mathbf{B}_1^{(u)}=(b_{ij}^{(u)})$ as a function of the environmental conditions characterizing the corresponding period of time.

Define $\Omega'$ as a set of possible states of environment and let $\mathbf{M}$ be a function mapping $\Omega'$ into the set of $L\times L$ stochastic matrices. Given a history of past environmental conditions $\omega=(\omega_1,\omega_2,\ldots)$ with $\omega_1\in\Omega',\omega_2\in\Omega',\ldots$ we put
\begin{equation}\label{stat}
  \mathbf{B}_1^{(u)}\equiv\mathbf{B}_1^{(u)}(\omega)=\mathbf{M}(\omega_u),\ u=1,2,\ldots.
\end{equation}
A simple choice of the state space $\Omega'=\{1,\ldots,K\}$ is a finite set with $K$ possible values for the random transition matrices $\mathbf{B}_1^{(u)}$.
Note that $K=1$ corresponds to the constant environment case. Two examples in Section \ref{ex} treat special cases with $K=2$ and $K=L$.

Our key assumption on the environmental history is that of stationarity
\[(\omega_1,\omega_2,\ldots)\stackrel{d}{=}(\omega_2,\omega_3,\ldots).\]
In this framework the fate of a single lineage is governed by the product of transition matrices
\begin{equation}\label{pu}
\mathbf{B}_1^{(1)}\cdots\mathbf{B}_1^{(u)}=\mathbf{M}(\omega_1)\cdots\mathbf{M}(\omega_u)=(p_{ij}^u)
\end{equation}
whose ergodic properties are well studied in  \cite{N}, \cite{C}, and \cite{O}.
An ergodic condition suitable for our purposes  is the following (see condition (D) on page 203 in \cite{C} and condition (a) on page 87 in \cite{N}):
\begin{align}\label{cN}
  &\mbox{for any $(i,j)$ and almost every realization of $\omega$ there exist a } \nonumber\\
&\mbox{$u=u_{ij}(\omega)$ and a $k=k_{ij}(\omega,u)$ such that the elements $p_{ik}^u$ and $p_{jk}^u$}\\
&\mbox{of the matrix \eqref{pu} are positive}.\nonumber
\end{align}

According to Theorem 6 in \cite{O} (see also Theorem 14 in \cite{N}), there exist {\it random} stationary probabilities $\gamma_i=\gamma_i(\omega),\ i=1,\ldots,L$  under condition \eqref{cN},  such that
\begin{equation*}
 \mathbf{B}_1^{(1)}\cdots\mathbf{B}_1^{(u)}\stackrel{d}{\to}
\mathbf{P}_1\equiv\left(
\begin{array}{ccc}
\gamma_1&\ldots&\gamma_L\\
\vdots&\ldots&\vdots\\
\gamma_1&\ldots&\gamma_L
\end{array}\right),\  u\to\infty
\end{equation*}
in distribution. Here the randomness of stationary probabilities for the single lineage position reflects environmental fluctuations. Next we state the main result of this paper allowing for dependence on $N$ in the sense of \eqref{N}: it is assumed that the backward migration probabilities have the form
\begin{equation}
  \label{N1}
\mathbf{B}_1^{(u)}(N)=\mathbf{B}_1^{(u)}+N^{-1}\mathbf{D}_1^{(u)}(N),
\end{equation}
where, as above, the matrices  $\mathbf{B}_1^{(u)}$ are genuine transition matrices, while the elements of the matrices $\mathbf{D}_1^{(u)}(N)$ are uniformly bounded in $u,N=1,2,\ldots$. Besides stationarity we will require the {\it mixing} property for the sequence of matrices  $\mathbf{B}_1^{(u)}$, meaning {\it asymptotic independence} between remote elements of the sequence (see Appendix A for technical details).

\begin{theorem}\label{main}
Consider a structured Wright-Fisher population with a random environment specified by the backward transition matrices $\mathbf{B}_1^{(u)}(N)$, $u=1,2,\ldots$ of the form \eqref{N1}. Assume that the sequence of matrices $\mathbf{B}_1^{(u)}$ is stationary and mixing. Under the condition \eqref{cN}, its ancestral process is approximated by the standard Kingman coalescent process
\begin{equation}
  \label{cc}
\{X([Nt/c_q]), t\ge0\}\to\{K(t), t\ge0\},\ N\to \infty,
\end{equation}
resulting in the coalescent effective population size formula \eqref{cq}.
\end{theorem}

In \eqref{c1}, \eqref{ce}, and \eqref{cc} convergence of stochastic processes is understood in the {\it Skorokhod sense} (which in this partcular setting is just a tiny improvement over convergence of finite-dimensional distributions). In these three coalescent approximation results the Skorokhod convergence follows from one-dimensional convergences like \eqref{11}, thanks to the Markov nature of the ancestral processes. The appropriate reference here is Theorem 2.12 on page 173 of \cite{EK}, called the Projection Theorem in \cite{NK}.

\section{Examples}\label{ex}
An important special case when the conditions of Theorem \ref{main} hold is that of random migration matrices $\mathbf{B}_1^{(u)}$ which are {\it independent and identically distributed} over $u=1,2,\ldots$. Then the path of a single lineage's is the trajectory of a Markov chain with random transition matrices, as
considered in \cite{T}. In the irreducible and aperiodic case, when
\begin{align}\label{irr}
  &\mbox{for any $(i,j)$ there is a $u=u_{ij}$ such that the element $p_{ij}^u$}\nonumber\\
&\mbox{of the random matrix \eqref{pu} satisfies }\mathbb{P}(p_{ij}^u>0)>0,
\end{align}
and
\begin{align}\label{st}
  \mathbb{P}(p_{ij}^u>0 \mbox{ for all } i)>0 \mbox{ for some $j$ and $u$},
\end{align}
the random vector of stationary probabilities $(\gamma_1,\ldots,\gamma_L)$ is strongly positive.

\begin{figure}
\centering
\includegraphics[width=10cm]{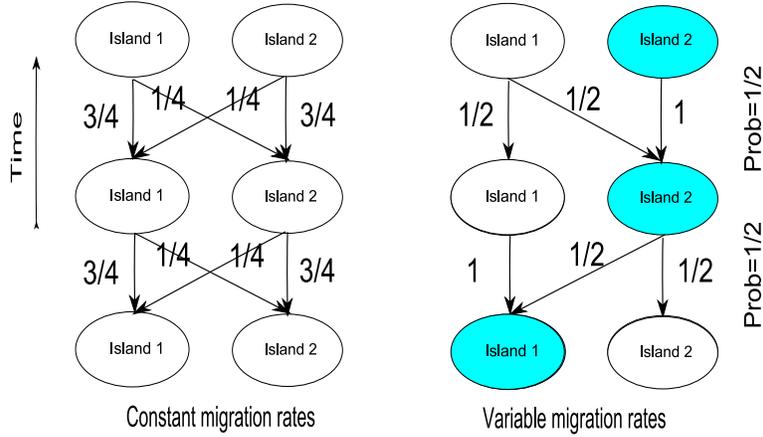}
\caption{The right panel presents a concrete example of a two-island model with variable migration. The random stationary probabilities for the backward mutation process have a uniform distribution. The left panel depicts the annealed version of the model with a fixed stationary distribution $\bar\gamma_1=\bar\gamma_2=0.5$.}\label{F2}
\end{figure}

This section contains two examples of population models with iid random environments which allow explicit calculations of products of transition matrices for migration processes. Our first example, illustrated by the right panel in Figure \ref{F2}, is a two-island ($L=2$) population model with an arbitrary $(a_1,a_2)$ satisfying \eqref{a}.
 Accordingly, the two sub-populations in a given generation consist of individuals labeled by numbers $1,\ldots,[Na_1]$ and $[Na_1]+1,\ldots,N$.

The defining one-step migration rules follow the next simple algorithm assuming just $K=2$ possible states of environment:
\begin{enumerate}
\item Toss a coin to decide which of the islands is favored environmentally,
\item If island 1 is favored, each of individuals  $1,\ldots,[N(a_1+\frac{a_2}{2})]$ chooses a parent uniformly at random from the previous generation individuals labeled  $1,\ldots,[Na_1]$, while each of individuals  $[N(a_1+\frac{a_2}{2})]+1,\ldots,N$ chooses a parent uniformly at random from the previous generation individuals, labeled  $[Na_1]+1,\ldots,N$,
\item If island 2 is favored, each of individuals  $[N(a_1/2)]+1,\ldots,N$ chooses a parent uniformly at random from the previous generation individuals labeled  $[Na_1]+1,\ldots,N$, while each of individuals  $1,\ldots,[N(a_1/2)]$ chooses a parent uniformly at random from the previous generation individuals labeled $1,\ldots,[Na_1]$.
\end{enumerate}
Notice that the proposed labelling of individuals within two sub-populations does not bring an unintended deterministic feature into the genetic drift dynamics, thanks to the underlying Wright-Fisher rules of genetic sampling.

The left panel of Figure \ref{F2} depicts the annealed version of the model with symmetric migration probabilities resulting in the stationary vector $(\bar\gamma_1,\bar\gamma_2)=(0.5,0.5)$. In view of \eqref{ca} this gives a benchmark factor
\begin{equation}\label{1}
c_a={1\over4}\left({1\over a_1}+{1\over a_2}\right)
\end{equation}
for the forthcoming effective population size formulas.

The beauty of this example lies in the full description of the products of independent matrices $\mathbf{B}_1^{(1)},\ldots,\mathbf{B}_1^{(u)}$ with the common distribution
\[\mathbb{P}\left(\mathbf{B}_1^{(1)}=\left(
\begin{array}{ll}
j2^{-1}&1-j2^{-1}\\
(j-1)2^{-1}&1-(j-1)2^{-1}
\end{array}\right)\right)=2^{-1},\ j=1,2.\]
The forward product $\mathbf{B}_1^{(1)}\cdots\mathbf{B}_1^{(u)}$ has a uniform distribution over $2^u$ matrices of the form
\begin{eqnarray*}
\left(
\begin{array}{ll}
j2^{-u}&1-j2^{-u}\\
(j-1)2^{-u}&1-(j-1)2^{-u}
\end{array}\right),\ j=1,\ldots,2^{u},
\end{eqnarray*}
which is verified by induction
\begin{eqnarray*}
&&\left({j\over 2^{u}},1-{j\over2^{u}}\right)
\left(
\begin{array}{cc}
1&0\\
1/2&1/2
\end{array}\right)=\left({j+2^{u}\over2^{u+1}},1-{j+2^{u}\over2^{u+1}}\right),\\
&&\left({j\over 2^{u}},1-{j\over2^{u}}\right)
\left(
\begin{array}{cc}
1/2&1/2\\
0&1
\end{array}\right)=\left({j\over2^{u+1}},1-{j\over2^{u+1}}\right).
\end{eqnarray*}

\begin{figure}
\centering
\includegraphics[width=10cm]{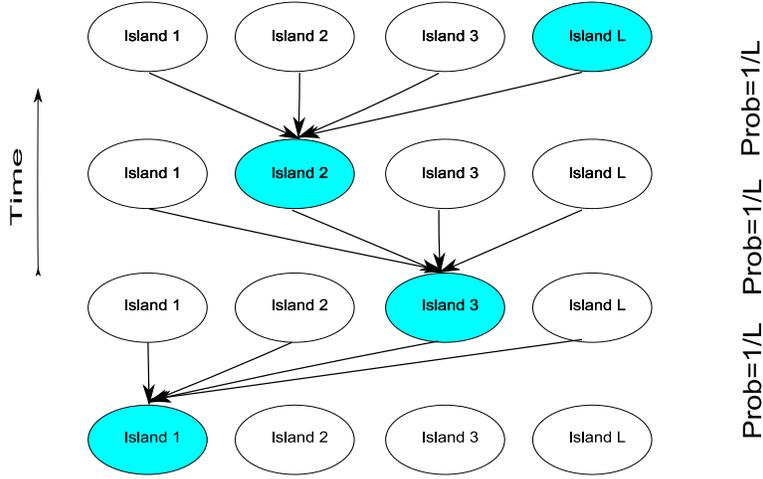}
\caption{A simple example of a multi-island model with variable migration.}\label{f3}
\end{figure}
The weak convergence of $\mathbf{B}_1^{(1)}\cdots\mathbf{B}_1^{(u)}$ as $u\to\infty$ (which is not an almost sure convergence) is made clear by the representation
\begin{eqnarray}\label{gex}
\mathbf{B}_1^{(1)}\cdots\mathbf{B}_1^{(u)}=\left(
\begin{array}{ll}
Z_u&1-Z_u\\
Z_u-2^{-u}&1-Z_u+2^{-u}
\end{array}\right),
\end{eqnarray}
where
\[Z_{u+1}=Z_u/2+\epsilon_u,\]
with iid $\epsilon_u$ taking values $0$ and $1/2$ with equal probabilities.

The reverse product has a similar representation
\begin{eqnarray*}
\mathbf{B}_1^{(u)}\cdots\mathbf{B}_1^{(1)}=\left(
\begin{array}{ll}
Z_u^*&1-Z_u^*\\
Z_u^*-2^{-u}&1-Z_u^*+2^{-u}
\end{array}\right)
\end{eqnarray*}
but with components
\[Z_{u+1}^*=Z_u^*-2^{-u}\epsilon_u\]
converging almost surely!
This remarkable phenomenon of different modes of convergence for different product orders of random matrices, well-known to mathematicians, might seem counterintuitive at first sight. The following simple observation may provide an illuminating parallel. Consider two sequences of random numbers $0.x_1x_2\ldots x_u$ and  $0.x_u\ldots x_2x_1$, where  $x_1,x_2\ldots$ are iid random digits. Clearly, the first sequence converges almost surely, and the second one only weakly, as $u\to\infty$. In both cases the limiting random number is uniformly distributed over the unit interval.

For this example it follows that the random stationary distribution vector $(\gamma_1,\gamma_2)$ has uniform components $\gamma_1\stackrel{d}{=}\gamma_2\sim U(0,1)$. Therefore, according to \eqref{cq} the corresponding factor for the quenched effective population size is given by
\begin{equation}\label{2}
 c_q^{(1)}={1\over3}\left({1\over a_1}+{1\over a_2}\right).
\end{equation}

Our second example is illustrated by Figure \ref{f3}. Now there is an arbitrary number $L$ of islands but migration rules are extremely simple. For each generation one island is chosen uniformly at random to be environmentally favored. Only the favored sub-population is giving offspring in the next generation as shown in the Figure \ref{f3}. In this case the stationary vector has a symmetric multivariate Bernoulli distribution $(\gamma_1,\ldots,\gamma_L)\sim$ Mn$(1,1/L,\ldots,1/L)$ resulting in the harmonic mean formula for the quenched effective population size
\begin{equation}\label{3}
 c_q^{(2)}= {1\over L}\sum_{k=1}^L{1\over a_k}.
\end{equation}
Viewed backwards in time, this example becomes a particular case of a much more general population model with variable population size considered in \cite{JS}.
\begin{figure}
\centering
\includegraphics[width=11cm,height=5cm]{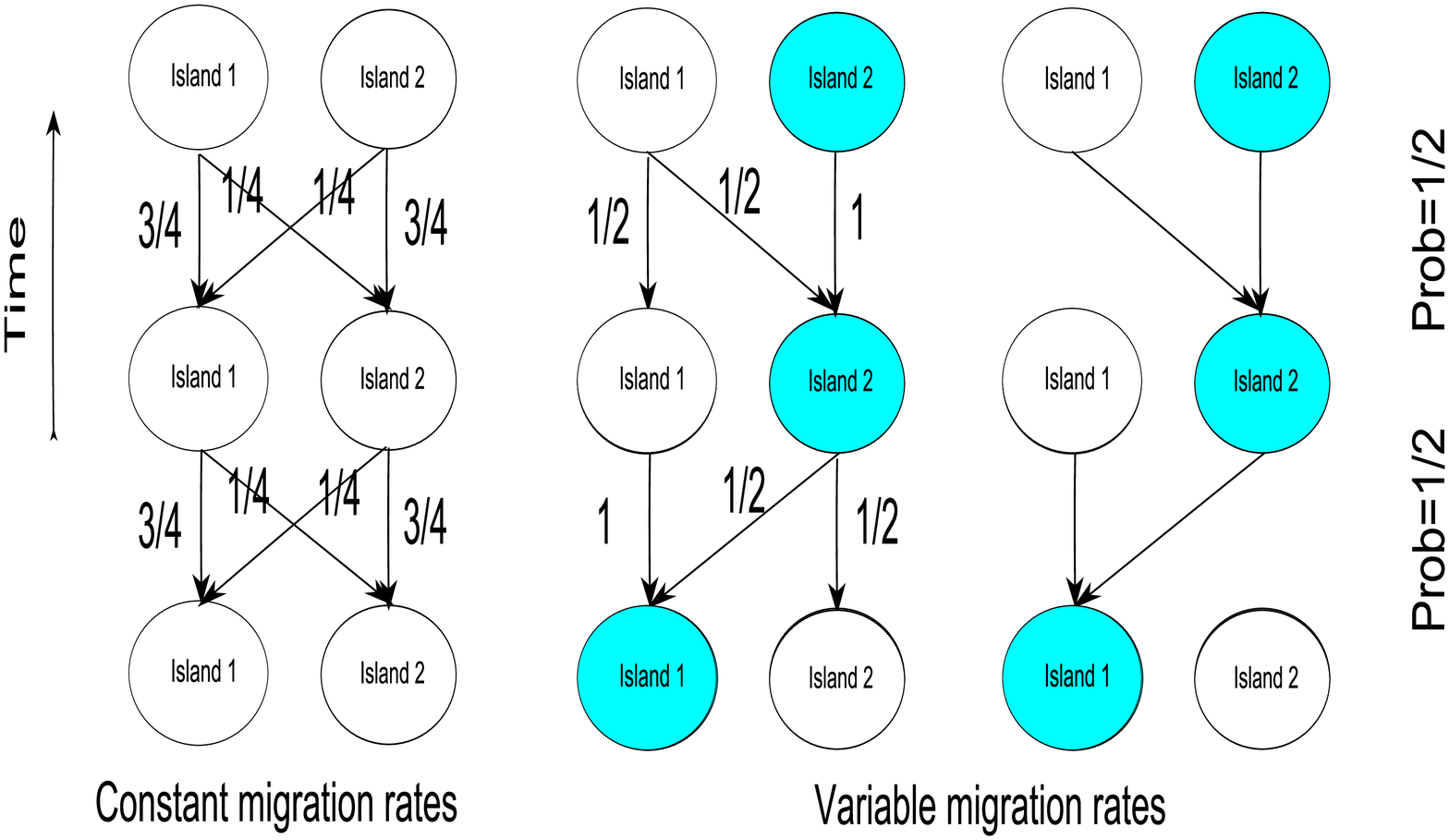}
\caption{Two population models with $L=2$ compared with the common annealed version.}\label{f4}
\end{figure}

Notice that for both examples conditions \eqref{irr} and \eqref{st} follow from
\[\mathbb{P}(b_{ij}^{(u)}>0 \mbox{ for all } i)>0,\ \  \mbox{ for all }j.\]
To summarize our examples we refer to Figure \ref{f4} which puts three sister models with two islands each together. Equations \eqref{1}-\eqref{3} yield the following correction formulas for two quenched $N_e$-factors as compared to the common annealed effective population size factor $c_a$:
\begin{eqnarray*}
&&  c^{(1)}_q={4\over3}c_a,\\
&& c^{(2)}_q=2c_a.
\end{eqnarray*}

\section{Annealed effective population size}\label{Sc}

As a prelude to the random environment case in Section \ref{pro}, a modified proof of \eqref{cf}, given in \cite{NK} will be outlined. In the current context, formula \eqref{cf} yields the annealed effective population size factor  \eqref{ca}, as explained in the introduction.

Throughout this section we assume constant environment and argue in terms of the configuration process of  $n$ lineages $\{\mathbf{X}(u)\}$, $\mathbf{X}(u)=(X_1(u),\ldots,X_L(u))$, where $X_i(u)$
is the number of lineages located on the $i$-th island at the $u$-th generation backward in time.
This is a Markov chain with the finite state space $S=\cup_{r=1}^n S_r$, where $S_r$ is the set of $r$-level
states $\mathbf{x}=(x_1,\ldots,x_L)$ with non-negative integer valued components satisfying $x_1+\cdots+x_L=r$. The number of elements in  $S_r$ is  $d_r={r+L-1\choose r}$.

Consider for a moment the backward migration process of $r$ lineages neglecting the possibility of coalescence. The corresponding transition matrix
$\mathbf{B}_r$ is of size $d_r\times d_r$. Since the Wright-Fisher
reproduction rule ensures that the paths of $r$ lineages are
independent, it is clear that the stationary distribution of the configuration
process on level $r$ is multinomial:
\begin{equation}
  \label{mn}
  \pi_r(\mathbf{x})={r\choose x_1,\ldots,x_L}\gamma_1^{x_1}\ldots\gamma_L^{x_L}, \ \mathbf{x}\in S_r,\ r=2,\ldots,n.
\end{equation}

For the transition matrix $\mathbf{\Pi}=(\Pi(\mathbf{x},\mathbf{y}))$ of the Markov chain $\{\mathbf{X}(u)\}$, the following counterpart of decomposition \eqref{k} is valid:
\begin{equation}
  \label{dec}
  \mathbf{\Pi}=\mathbf{B}(\mathbf{I}+N^{-1}\mathbf{C})+\mathbf{o}(N^{-1}),
\end{equation}
where $\mathbf{B}={\rm diag}(\mathbf{B}_1,\ldots,\mathbf{B}_n)$ is
the block diagonal matrix with the transition probabilities caused by pure migration (coalescence prohibited), while the matrix $\mathbf{C}$ gives the coalescence rates for various geographical configurations of sampled ancestral lines,
\begin{equation}
  \label{C}
 \mathbf{C}=
\left(
\begin{array}{lllllll}
\mathbf{O}_{11}&\mathbf{O}_{12}&\mathbf{O}_{13}&\ldots&\mathbf{O}_{1,n-2}&\mathbf{O}_{1,n-1}&\mathbf{O}_{1n}\\
\mathbf{C}_{21}&-\mathbf{C}_2&\mathbf{O}_{23}&\ldots&\mathbf{O}_{2,n-2}&\mathbf{O}_{2,n-1}&\mathbf{O}_{2n}\\
\mathbf{O}_{3,1}&\mathbf{C}_{3,2}&-\mathbf{C}_{3}&\ldots&\mathbf{O}_{3,n-2}&\mathbf{O}_{3,n-1}&\mathbf{O}_{3,n}\\
\vdots&\vdots&\vdots&\ldots&\vdots&\vdots&\vdots\\
\mathbf{O}_{n1}&\mathbf{O}_{n2}&\mathbf{O}_{n3}&\ldots&\mathbf{O}_{n,n-2}&\mathbf{C}_{n,n-1}&-\mathbf{C}_{n}
\end{array}\right).
\end{equation}
Here matrices $\mathbf{O}_{ij}$ have dimensions $d_i\times d_j$ and their elements are all zero. The  blocks $\mathbf{\mathbf{C}}_{r}$ on the main diagonal of $\mathbf{C}$ are diagonal matrices themselves,
$$\mathbf{\mathbf{C}}_{r}={\rm diag}(C(\mathbf{x}), \mathbf{x}\in S_r),\ C(\mathbf{x})=\sum_{k=1}^L{1\over a_k}{x_k\choose2}.$$
The $d_r\times d_{r-1}$ blocks $\mathbf{\mathbf{C}}_{r,r-1}$, constituting the ``second diagonal'' of $\mathbf{C}$, have elements ${1\over a_k}{x_k\choose2}$ at positions  $(\mathbf{x},\mathbf{x}-\mathbf{e}_k)$ with
\begin{align*}
\mathbf{x}&=(x_1,\ldots,x_L)\in S_r,\\
\mathbf{x}-\mathbf{e}_k&=(x_1,\ldots,x_{k-1},x_k-1,x_{k+1},\ldots,x_L)\in S_{r-1},
\end{align*}
and zero elements elsewhere.

In particular, if $L=2$, then $d_r=r+1$ counts two dimensional configurations, which we will order in the following way: $(r,0),(r-1,1),\ldots,(1,r-1),(0,r)$. The non-zero blocks of the matrix $\mathbf{C}$ are of two kinds: the $(r+1)\times(r+1)$ matrices
\begin{equation*}
 \mathbf{C}_r=
\left(
\begin{array}{cccccccc}
{r\choose2}{1\over a_1}&0&0&\ldots&0&\ldots&0&0\\
\vdots&\vdots&\vdots&\ldots&\vdots&\ldots&\vdots&\vdots\\
0&0&0&\ldots&{r-k\choose2}{1\over a_1}+{k\choose2}{1\over a_2}&\ldots&0&0\\
\vdots&\vdots&\vdots&\ldots&\vdots&\ldots&\vdots&\vdots\\
0&0&0&\ldots&0&\ldots&0&{r\choose2}{1\over a_2}
\end{array}\right)
\end{equation*}
and $(r+1)\times r$ matrices
\begin{equation*}
 \mathbf{C}_{r,r-1}=
\left(
\begin{array}{cccccccc}
{r\choose2}{1\over a_1}&0&0&\ldots&0&0&\ldots&0\\
0&{r-1\choose2}{1\over a_1}&0&\ldots&0&0&\ldots&0\\
0&{1\over a_2}&{r-2\choose2}{1\over a_1}&\ldots&0&0&\ldots&0\\
\vdots&\vdots&\vdots&\ldots&\vdots&\vdots&\ldots&\vdots\\
0&0&0&\ldots&{k\choose2}{1\over a_2}&{r-k\choose2}{1\over a_1}&\ldots&0\\
\vdots&\vdots&\vdots&\ldots&\vdots&\vdots&\ldots&\vdots\\
0&0&0&\ldots&0&0&\ldots&{r\choose2}{1\over a_2}
\end{array}\right).
\end{equation*}

Put $X(u)=X_1(u)+\cdots+X_L(u)$. What we are really interested in is not the Markov chain $\{\mathbf{X}(u)\}$ itself, but rather its collapsed version $\{X(u)\}$ focusing on the total number of lineages and disregarding the frequently changing geographical locations of the sampled lineages. Clearly, the total number of lineages $X(u)$ is not generally a Markov process.
%except for the "dummy islands" case discussed in the end of this section. 
Given a matrix $\mathbf{R}=(R(\mathbf{x},\mathbf{y}))$ of the same dimension $(d_1+\cdots+d_n)\times(d_1+\cdots+d_n)$ as the matrix $\mathbf{\Pi}$, we write $\mathbf{R}^{\downarrow}_{\mathbf{x}_1,\ldots,\mathbf{x}_n}$ to denote its collapsed version of size $n\times n$ with elements $\sum_{\mathbf{y}\in S_j}R(\mathbf{x}_i,\mathbf{y})$
depending on a specified set of elements $\mathbf{x}_i\in S_i$, $1\le i,j\le n$.
In this notation, the desired convergence to the Kingman coalescent \eqref{ce} is equivalent to the claim that for any given vector $(\mathbf{x}_1,\ldots,\mathbf{x}_n)$
\begin{equation}
  \label{ki}
(\mathbf{\Pi}^{Nt})^{\downarrow}_{\mathbf{x}_1,\ldots,\mathbf{x}_n}\to e^{c_ft\mathbf{Q}},\ N\to\infty.
\end{equation}

%Convergence \eqref{ki} 
This follows from M\"ohle's lemma \cite{M} which in view of \eqref{dec}  gives
\begin{equation}
  \label{m1}
  \mathbf{\Pi}^{Nt}\to\mathbf{P}-\mathbf{I}+e^{t\mathbf{P}\mathbf{C}\mathbf{P}},\ N\to\infty,
\end{equation}
where $\mathbf{P}={\rm diag}(\mathbf{P}_1,\ldots,\mathbf{P}_n)$ is a block diagonal matrix with the block $\mathbf{P}_r$ having $d_r$ equal rows $(\pi_r(\mathbf{x}), \mathbf{x}\in S_r)$. To reconcile \eqref{ki} and \eqref{m1} we suggest using a representation
\begin{equation}\label{star}
 \mathbf{P}-\mathbf{I}+ e^{t\mathbf{P}\mathbf{C}\mathbf{P}}=e^{c_ft\mathbf{Q}}*\mathbf{P},
\end{equation}
where a special matrix product  $\mathbf{G}*\mathbf{P}$ is defined for an arbitrary $n\times n$ matrix $\mathbf{G}=(g_{ij})$ as a block matrix with blocks  $g_{ij}\mathbf{P}_{i,j}$, where  $\mathbf{P}_{i,j}$ has  $d_i$ rows, each equal to a row in $\mathbf{P}_{j}$:
\begin{equation*}
 \mathbf{G}*\mathbf{P}=
\left(
\begin{array}{c|c|c|c}
g_{11}\pi_1(\mathbf{x}), \mathbf{x}\in S_1&g_{12}\pi_2(\mathbf{x}), \mathbf{x}\in S_2&\ldots&g_{1n}\pi_n(\mathbf{x}), \mathbf{x}\in S_n\\
\vdots&\vdots&\ldots&\vdots\\
g_{11}\pi_1(\mathbf{x}), \mathbf{x}\in S_1&g_{12}\pi_2(\mathbf{x}), \mathbf{x}\in S_2&\ldots&g_{1n}\pi_n(\mathbf{x}), \mathbf{x}\in S_n\\\hline
g_{21}\pi_1(\mathbf{x}), \mathbf{x}\in S_1&g_{22}\pi_2(\mathbf{x}), \mathbf{x}\in S_2&\ldots&g_{2n}\pi_n(\mathbf{x}), \mathbf{x}\in S_n\\
\vdots&\vdots&\ldots&\vdots\\
g_{21}\pi_1(\mathbf{x}), \mathbf{x}\in S_1&g_{22}\pi_2(\mathbf{x}), \mathbf{x}\in S_2&\ldots&g_{2n}\pi_n(\mathbf{x}), \mathbf{x}\in S_n\\\hline
\vdots&\vdots&\ldots&\vdots\\\hline
g_{n1}\pi_1(\mathbf{x}), \mathbf{x}\in S_1&g_{n2}\pi_2(\mathbf{x}), \mathbf{x}\in S_2&\ldots&g_{nn}\pi_n(\mathbf{x}), \mathbf{x}\in S_n\\
\vdots&\vdots&\ldots&\vdots\\
g_{n1}\pi_1(\mathbf{x}), \mathbf{x}\in S_1&g_{n2}\pi_2(\mathbf{x}), \mathbf{x}\in S_2&\ldots&g_{nn}\pi_n(\mathbf{x}), \mathbf{x}\in S_n
\end{array}\right).
\end{equation*}
Clearly,
$$(\mathbf{G}*\mathbf{P})^{\downarrow}_{\mathbf{x}_1,\ldots,\mathbf{x}_n}=\mathbf{G},$$
irrespective of the choice of  $(\mathbf{x}_1,\ldots,\mathbf{x}_n)$.

To verify \eqref{star}, notice that the product $\mathbf{P}\mathbf{C}\mathbf{P}$ has the same structure as the matrix  $\mathbf{C}$ with
blocks $(-\mathbf{P}_r\mathbf{C}_{r}\mathbf{P}_r)$ on the main diagonal and blocks
$\mathbf{P}_r\mathbf{C}_{r,r-1}\mathbf{P}_{r-1}$ on the second diagonal. This observation together with
\begin{eqnarray*}
  \sum_{\mathbf{y}\in S_r}\pi_r(\mathbf{y})C(\mathbf{y})&=&\sum_{\mathbf{y}\in S_r}{r\choose y_1,\ldots,y_L}\gamma_1^{y_1}\cdots\gamma_L^{y_L}\sum_{k=1}^L{1\over a_k}{y_k\choose2}\\
&=&\sum_{k=1}^L{1\over 2a_k}\sum_{\mathbf{y}\in S_r}y_k(y_k-1){r\choose y_1,\ldots,y_L}\gamma_1^{y_1}\cdots\gamma_L^{y_L}\\
&=&\sum_{k=1}^L{\gamma_k^2\over 2a_k}{\partial^2\over\partial s_k^2}(s_1+\cdots+s_L)^r|_{s_1=\gamma_1,\ldots,s_L=\gamma_L}={r\choose2}c_f
\end{eqnarray*}
implies that
\begin{equation}\label{pcp}
\mathbf{P}\mathbf{C}\mathbf{P}=c_f\mathbf{Q}*\mathbf{P}.
\end{equation}
It remains to observe that $(\mathbf{Q}*\mathbf{P})^k=\mathbf{Q}^k*\mathbf{P}$.

\section{Proof of Theorem \ref{main}}\label{pro}

Without loss of generality the sequence of environmental states can be viewed as a doubly infinite stationary sequence
\[\omega=(\ldots,\omega_{-2},\omega_{-1},\omega_{0},\omega_{1},\omega_{2},\ldots).\]
According to Theorem 6 in \cite{C} (see also Theorem 14 in \cite{N}), condition \eqref{cN} guarantees that the matrix product in the reversed order converges almost surely as $u\to\infty$:
\begin{equation}
  \label{ergr}
 \mathbf{B}_1^{(-u)}\cdots\mathbf{B}_1^{(-1)}\mathbf{B}_1^{(0)}\mathbf{B}_1^{(1)}\cdots\mathbf{B}_1^{(j)}\to%\Gamma=
\mathbf{P}_1^{(j)}\equiv\left(
\begin{array}{ccc}
\gamma_1^{(j)}&\ldots&\gamma_L^{(j)}\\
\vdots&\ldots&\vdots\\
\gamma_1^{(j)}&\ldots&\gamma_L^{(j)}
\end{array}\right).
\end{equation}
Importantly, the vectors $(\gamma_1^{(j)},\ldots,\gamma_L^{(j)})\stackrel{d}{=}(\gamma_1,\ldots,\gamma_L)$ satisfy a recursive relation
\begin{equation}\label{rec}
(\gamma_1^{(j)},\ldots,\gamma_L^{(j)})\mathbf{B}_1^{(j+1)}=(\gamma_1^{(j+1)},\ldots,\gamma_L^{(j+1)}).
\end{equation}

Let $\mathbf{B}^{(j)}={\rm diag}(\mathbf{B}_1^{(j)},\ldots,\mathbf{B}_n^{(j)})$, $j=1,2,\ldots$ be the block diagonal matrices characterizing configurations of non-coalescing lineages. We have weak convergence of random matrices
\begin{equation}\label{bi}
  \mathbf{B}^{(1)}\cdots\mathbf{B}^{(u)}\stackrel{d}{\to}\mathbf{P},\  u\to\infty,
\end{equation}
where $\mathbf{P}$ is defined by \eqref{mn} in terms of the (now random) vector $(\gamma_1,\ldots,\gamma_L)$ exactly as in Section \ref{Sc}. On other hand, we can rely on the a.s. convergence
\begin{equation}\label{bip}
\mathbf{B}^{(-u)}\cdots\mathbf{B}^{(-1)}\mathbf{B}^{(0)}\mathbf{B}^{(1)}\cdots\mathbf{B}^{(j)}\stackrel{a.s.}{\to}\mathbf{P}^{(j)},\  u\to\infty,\ j\ge1,
\end{equation}
where $\mathbf{P}^{(j)}\stackrel{d}{=}\mathbf{P}$ are all defined on the same probability space using vectors $(\gamma_1^{(j)},\ldots,\gamma_L^{(j)})$ given by \eqref{ergr}. Observe that
since the rows of matrix $\mathbf{P}^{(j)}$ are identical, we have
\begin{equation}\label{pipi}
\mathbf{P}^{(i)}\mathbf{P}^{(j)}=\mathbf{P}^{(j)},
\end{equation}
for any pair $(i,j)$, and moreover, due to \eqref{rec}, we have a recursion
\begin{equation}\label{pip}
\mathbf{P}^{(j)}\mathbf{B}^{(j+1)}=\mathbf{P}^{(j+1)}.
\end{equation}

 The proof of Theorem \ref{main} extends the approach outlined in the previous section and establishes the following almost sure convergence of random transition probabilities for the configuration process $\mathbf{X}(u)=(X_1(u),\ldots,X_L(u))$:
\begin{equation}
  \label{kik}
\|\mathbf{\Pi}^{(1)}_N\cdots\mathbf{\Pi}^{([Nt])}_N-e^{c_qt\mathbf{Q}}*\mathbf{P}^{([Nt])}\|\stackrel{a.s.}{\to} 0,\ N\to\infty,
\end{equation}
where the norm of a matrix $\mathbf{G}=(g_{ij})$ is defined as $\left\Vert \mathbf{G}\right\Vert=\max_{i}\sum_{j}|g_{ij}|$. As we show in Appendix B, this follows from the next two key observations:
\begin{equation}
  \label{kiki}
\|\mathbf{B}^{(1)}\cdots\mathbf{B}^{(u)}-\mathbf{P}^{(u)}\|\stackrel{a.s.}{\to} 0,\ u\to\infty,
\end{equation}
where $\mathbf{P}^{(u)}$ are defined by \eqref{bip}, and (see Appendix A)
\begin{equation}\label{er}
{1\over u}\sum_{j=1}^{u}c^{(j)}\stackrel{a.s.}{\to} \mathbb{E}\left(\sum_{k=1}^L{1\over a_k}\gamma_k^2\right)=c_q,\ u\to\infty,
\end{equation}
where
\begin{equation}\label{cj}
c^{(j)}=\sum_{k=1}^L{1\over a_k}(\gamma^{(j)}_k)^2.
\end{equation}
Observe that for the first example in Section \ref{ex} the almost sure convergence \eqref{kiki} holds with $(\gamma_1^{(u)},\gamma_2^{(u)})=(Z_u,1-Z_u)$, due to \eqref{gex}.

To prove \eqref{kiki} observe that due to \eqref{pip} the random sequence $\Delta_u=\|\mathbf{B}^{(1)}\cdots\mathbf{B}^{(u)}-\mathbf{P}^{(u)}\|$ is monotone:
\[\Delta_{u+1}\stackrel{\eqref{pip}}{=}\|(\mathbf{B}^{(1)}\cdots\mathbf{B}^{(u)}-\mathbf{P}^{(u)})\mathbf{B}^{(u+1)}\|\le\|\mathbf{B}^{(1)}\cdots\mathbf{B}^{(u)}-\mathbf{P}^{(u)}\|=\Delta_u.\]
It remains to note the convergence in probability $\Delta_u\stackrel{P}{\to} 0$,
which follows from \eqref{bi} and the representation $\mathbf{P}^{(u)}=\mathbf{P}^{(0)}\mathbf{B}^{(1)}\cdots\mathbf{B}^{(u)}$ implied by \eqref{bip}.

\section{Appendix A: proof of \eqref{er}}

To verify  \eqref{er} it is enough to check that for any $k=1,\ldots,L$
\[{1\over u}\sum_{j=1}^{u}\left(\gamma_k^{(j)}\right)^2\stackrel{a.s.}{\to} \mathbb{E}\left(\gamma_k^2\right),\ u\to\infty.\]
According to the ergodic theorem discussed in Chapter 6.4 of \cite{D}, this would follow if we show that the stationary sequence $\{\gamma_k^{(u)}\}_{u=-\infty}^\infty$ posesses the mixing property (remote elements of the sequence are asymptotically independent):
\begin{equation}\label{er1}
\mathbb{P}(\gamma_k^{(0)}\le x;\gamma_k^{(-u)}\le y)\to \mathbb{P}(\gamma_k^{(0)}\le x)\mathbb{P}(\gamma_k^{(0)}\le y),\ u\to\infty.\
\end{equation}
whatever are $x\in[0,1]$ and $y\in[0,1]$. As we show next, relation \eqref{er1} follows from the representation
\begin{equation}\label{rep}
(\gamma_1^{(-u)},\ldots,\gamma_L^{(-u)})\mathbf{B}_1^{(-u+1)}\cdots\mathbf{B}_1^{(0)}=(\gamma_1^{(0)},\ldots,\gamma_L^{(0)}),
\end{equation}
see \eqref{rec}, and
the assumed  mixing property for the sequence of matrices $\mathbf{B}_1^{(u)}$. The latter says that any two events separated by a large number $u$ of units of time
\begin{align*}
&A\in\sigma\{\mathbf{B}_1^{(0)},\mathbf{B}_1^{(1)},\ldots\},\\
&A_u\in\sigma\{\mathbf{B}_1^{(-u)},\mathbf{B}_1^{(-u-1)},\ldots\},
\end{align*}
with $\mathbb{P}(A)=p_1$ and $\mathbb{P}(B_u)=p_2$, are asymptotically independent:
\[\mathbb{P}(A\cap B_u)\to p_1p_2,\ u\to\infty.\]
(Here $\sigma\{\xi_1,\xi_2,\ldots\}$ stands for the sigma-algebra of events generated by the random variables $\xi_1,\xi_2,\ldots$.)

Define $\alpha_k^{(u)}\le \beta_k^{(u)}$ as the minimal and maximal elements in the $k-$th column of the matrix product $\mathbf{B}_1^{(-u+1)}\cdots\mathbf{B}_1^{(0)}$. It is easily verified that $\alpha_k^{(u)}$ increases while  $\beta_k^{(u)}$  decreases with $u$ since each realization of $\mathbf{B}_1^{(u)}$ is a stochastic matrix (every row is a non-negative vector with components summing to 1). Clearly, for any natural $v$,
\[\mathbb{P}(\gamma_k^{(0)}\le x;\gamma_k^{(-u)}\le y)\le \mathbb{P}(\alpha_k^{(v)}\le x;\gamma_k^{(-u)}\le y),\]
which due to stationarity of $\{\mathbf{B}_1^{(u)}\}_{u=-\infty}^\infty$ and its mixing property implies
\[\limsup_{u\to\infty}\mathbb{P}(\gamma_k^{(0)}\le x;\gamma_k^{(-u)}\le y)\le \mathbb{P}(\alpha_k^{(v)}\le x)\mathbb{P}(\gamma_k^{(0)}\le y).\]
As we already know, under condition  \eqref{cN},  $\alpha_k^{(v)}\nearrow\gamma_k^{(0)}$ as $v\to\infty$. Thus $\{\alpha_k^{(v)}\le x\}\searrow\{\gamma_k^{(0)}\le x\}$ and it follows that
\[\limsup_{u\to\infty}\mathbb{P}(\gamma_k^{(0)}\le x;\gamma_k^{(-u)}\le y)\le \mathbb{P}(\gamma_k^{(0)}\le x)\mathbb{P}(\gamma_k^{(0)}\le y).\]
A similar reasoning in terms of $\beta_k^{(v)}$ gives the lower bound
\[\liminf_{u\to\infty}\mathbb{P}(\gamma_k^{(0)}\le x;\gamma_k^{(-u)}\le y)\ge \mathbb{P}(\gamma_k^{(0)}\le x)\mathbb{P}(\gamma_k^{(0)}\le y)\]
finishing the proof of  \eqref{er1} and therefore of \eqref{er}.

\section{Appendix B: proof of \eqref{kik}}

Appropriately modifying the notation from Section \ref{Sc}, we set the starting point of our proof of \eqref{kik} in a form similar to \eqref{dec}:
\begin{equation*}
  \label{decr}
  \mathbf{\Pi}^{(j)}_N=\left(\mathbf{B}^{(j)}+N^{-1}\mathbf{D}^{(j)}(N)\right)(\mathbf{I}+N^{-1}\mathbf{C})+\mathbf{o}(N^{-1}).
\end{equation*}
Here elements of the matrices $\mathbf{D}^{(j)}(N)$ are uniformly bounded in $u$ and $N$, with all the rows of $\mathbf{D}^{(j)}(N)$ having zero sums. Thus,
\begin{eqnarray*}
\mathbf{\Pi} _{N}^{(j)}=\mathbf{B}^{(j)}+N^{-1}\mathbf{H}_{N}^{(j)},
\label{reprmatrix}
\end{eqnarray*}%
where
\begin{equation*}
\mathbf{H}_{N}^{(j)}=\mathbf{D}^{(j)}(N)+\mathbf{B}^{(j)}\mathbf{C}+\mathbf{o}(N^{-1}).
\end{equation*}
In view of \eqref{kiki} it is a straightforward exercise to modify the proof of the first part of Lemma 2.1 in \cite{M} to obtain
\begin{equation*}
  \label{kim}
\|\mathbf{\Pi}^{(1)}_N\cdots\mathbf{\Pi}^{([Nt])}_N-\prod_{j=1}^{[Nt]}(\mathbf{P}^{(j)}+N^{-1}\mathbf{H}_{N}^{(j)})\|\stackrel{a.s.}{\to} 0,\ N\to\infty.
\end{equation*}
Since the rows of $\mathbf{P}^{(j)}$ are identical and the rows of $\mathbf{D}^{(j)}_N$ sum to zero, we get
\[\mathbf{H}_{N}^{(j-1)}\mathbf{P}^{(j)}=\mathbf{B}^{(j-1)}\mathbf{C}\mathbf{P}^{(j)}+\mathbf{o}(N^{-1}).\]
Using \eqref{pipi} and \eqref{pip} we obtain
\begin{equation*}
  \label{kimi}
\|\prod_{j=1}^{[Nt]}(\mathbf{P}^{(j)}+N^{-1}\mathbf{H}_{N}^{(j)})-\prod_{j=1}^{[Nt]}(\mathbf{P}^{(j)}+N^{-1}\mathbf{B}^{(j)}\mathbf{C})\|\stackrel{a.s.}{\to} 0,\ N\to\infty,
\end{equation*}
and also
\begin{equation*}
  \label{kimim}
\|\prod_{j=1}^{[Nt]}(\mathbf{P}^{(j)}+N^{-1}\mathbf{B}^{(j)}\mathbf{C})-\prod_{j=1}^{[Nt]}\mathbf{P}^{(j)}(\mathbf{I}+N^{-1}\mathbf{C})\|\stackrel{a.s.}{\to} 0,\ N\to\infty.
\end{equation*}
Further, it is not difficult to check that for any $i$ and $j$
\begin{align*}
\left(\mathbf{Q}\ast \mathbf{P}^{(i)}\right)\left(\mathbf{Q}\ast \mathbf{P}^{(j)}\right)
&=\mathbf{Q}^{2}\ast \mathbf{P}^{(j)},\\
\left(\mathbf{Q}\ast \mathbf{P}^{(i)}\right)\mathbf{P}^{(j)}&=\mathbf{Q}\ast \mathbf{P}^{(j)},\\
\mathbf{P}^{(j)}\mathbf{C}\mathbf{P}^{(j)}&=c^{(j)}\mathbf{Q}\ast \mathbf{P}^{(j)}
\end{align*}%
(see \eqref{pcp} and \eqref{cj} for an explanation of the last equality).

For $m_{1},\ldots,m_{k}\in \mathbb{N}_{0}$ set $M_{0}=0,$ $M_{j}=m_{1}+\cdots+m_{j}$. Then
\begin{align*}
\prod_{j=1}^{[Nt]-1}&\mathbf{P}^{(j)}(\mathbf{I}+N^{-1}\mathbf{C})\mathbf{P}^{([Nt])}\\
&=\sum\limits_{k=1}^{[Nt]-1}\frac{1}{N^{k-1}}\sum_{\substack{ %
m_{1},...,m_{k}\in \mathbb{N}_0 \\ M_{k}=[Nt]-1}}%
\mathbf{P}^{(M_{1})}\mathbf{C}\mathbf{P}^{(M_{2})}\mathbf{C}\cdots \mathbf{P}^{(M_{k-1})}\mathbf{C}\mathbf{P}^{([Nt])} \\
&=\sum\limits_{k=1}^{[Nt]-1}\frac{1}{N^{k-1}}\sum_{\substack{ %
m_{1},...,m_{k}\in \mathbb{N}_0 \\  M_{k}=[Nt]-1}}%
\left(\prod_{j=1}^{k-1}\mathbf{P}^{(M_{j})}\mathbf{C}\mathbf{P}^{(M_{j})}\right)\mathbf{P}^{([Nt])} \\
&=\sum\limits_{k=1}^{[Nt]-1}\frac{1}{N^{k-1}}\sum_{\substack{ %
m_{1},...,m_{k}\in \mathbb{N}_0 \\  M_{k}=[Nt]-1}}%
\left(\prod_{j=1}^{k-1}c^{(M_{j})}\mathbf{Q}*\mathbf{P}^{(M_{j})}\right)\mathbf{P}^{([Nt])} \\
&=\sum\limits_{k=1}^{[Nt]-1}\sum_{\substack{ %
m_{1},...,m_{k}\in \mathbb{N}_0 \\  M_{k}=[Nt]-1}}%
\left(\prod_{j=1}^{k-1}{c^{(M_{j})}\over N}\mathbf{Q}\right)*\mathbf{P}^{([Nt])} \\
&=\prod_{j=1}^{[Nt]-1}\left(\mathbf{I}+{c^{(j)}\over N}\mathbf{Q}\right)
\ast \mathbf{P}^{([Nt])}.
\end{align*}%
To derive \eqref{kik} from the previous chain of relations it remains to observe that
\begin{equation*}
\left\|\prod_{j=1}^{[Nt]-1}\left(\mathbf{I}+c^{(j)}N^{-1}\mathbf{Q}\right)-\exp\left\{\left(N^{-1}\sum_{j=1}^{[Nt]}c^{(j)}\right)\mathbf{Q}\right\}\right\|\stackrel{a.s.}{\to}0,\ N\rightarrow \infty
\end{equation*}%
and apply \eqref{er}.


\begin{thebibliography}{99}

\bibitem{C} {\sc Cogburn, R.} (1986) On products of random matrices. Contemp. Math. 50, 199--213.
\bibitem{crowkimura}
         {\sc Crow, J.F. and Kimura, M.} (1970)
         An Introduction to Population Genetics Theory.
         Harper and Row, New York.
      \bibitem{D} {\sc Durrett, R.} (2004) Probability: Theory and Examples. Third Edition, Duxbury Press.
\bibitem{EK} {\sc Ethier, S.N. and Kurtz, T.G.} (1986) Markov
  Processes: Characterisation and Convergence. Wiley, New York.
    \bibitem{E79}
         {\sc Ewens, W.J.} (1979) Mathematical Population Genetics.
         Springer, Berlin.
\bibitem{E82} {\sc Ewens, W.J.} (1982) On the concept of effective size. Theor. Pop. Biol. 21, 373–-378.
\bibitem{E89} {\sc Ewens, W.J.} (1989) The effective population size in the presence of catastrophes, pp. 9–-25 in Mathematical Evolutionary Theory, edited by M. W. Feldman. Princeton University Press, Princeton, NJ.
\bibitem{E} {\sc Ewens, W.J.} (2004) Mathematical Population Genetics (2nd Edition). Springer-Verlag, New York.
\bibitem{JS} {\sc Jagers, P. and Sagitov, S.} (2004) Convergence to the coalescent in populations of substantially varying size. J. Appl. Prob. 41, 368--378.
\bibitem{KK} {\sc Kaj, I. and Krone, S.M.} (2003) The coalescent
  process in a population with stochastically varying size. J.
Appl. Prob. 40, 33--48.
\bibitem{K1} {\sc Kingman, J.F.C.} (1982) On the genealogy of large
  populations.  J. Appl. Prob. 19A, 27--43.
\bibitem{M} {\sc M\"ohle, M.} (1998) A convergence theorem for Markov chains arising in population genetics and the coalescent with selfing.  Adv. Appl. Prob. 30, 493--512.
      \bibitem{M2}
         {\sc M\"ohle, M.} (2001)
         Forward and backward diffusion approximations for haploid
         exchangeable population models.
         Stoch. Process. Appl. 95, 133--149.
\bibitem{Na} {\sc Nagylaki, T.} (1980) The strong-migration limit in geographically structured populations. J. Math. Biol. 9, 101--114.

\bibitem{N} {\sc Nawrotzki, K.} (1981-1982) Discrete open systems or Markov chains in a random environment. I,II. J. Inform. Process. Cybernet. 17, 569-599;  18, 83--98.

\bibitem{NT} {\sc Nei, M. and N. Takahata} (1993) Effective population size, genetic diversity, and coalescence time in subdivided populations. J. Mol. Evol. 37, 240–-244.
\bibitem{NK} {\sc Nordborg, M. and Krone, S.} (2002) Separation of
  time scales and convergence to the coalescent in structured
  populations. In Modern Developments in Theoretical
  Population Genetics, pp. 194--232, M. Slatkin and M.  Veuille, editors. Oxford
  University Press.
\bibitem{No} {\sc Notohara, M.} (1993) The strong-migration limit for the genealogical process in geographically
structured populations. J. Math. Biol. 31, 115--122.
\bibitem{O} {\sc Orey, S.} (1991) Markov chains with stochastically stationary transition probabilities. Ann. Prob. 19, 907--928.
\bibitem{SJ} {\sc Sagitov, S. and Jagers, P.} (2005) The coalescent effective size of age-structured populations. Ann. Appl. Prob. 15, 1778--1797.
\bibitem{SKK}{\sc Sj\"odin, P., Kaj,  I., Krone, S., Lascoux, M., and Nordborg, M.} (2004)
On the meaning and existence of an effective population size. Genetics 169, 1061--1070.
\bibitem{T} {\sc Takahashi, Y.} (1969)
Markov chains with random transition matrices.
Kodai Math. Sem. Rep. 21, 426--447.
\bibitem{W98} {\sc Wakeley, J.} (1998) Segregating sites in Wright's island model. Theor. Pop. Biol. 53, 166--174.
      \bibitem{WS}
         {\sc Wakeley, J. and Sargsyan, O.} (2009)
         Extensions of the coalescent effective population size. Genetics 181, 341--345.
  \end{thebibliography}
\end{document}